\documentclass[11pt]{amsart}
\usepackage{amsfonts}
\usepackage{amsmath}
\usepackage{calc}
\usepackage{graphicx}
\usepackage{amsmath}
\usepackage{amsfonts}
 \usepackage{amssymb}
 \usepackage{verbatim} 
\usepackage{amsthm} 
\usepackage[colorlinks=true, pdfstartview=FitV, linkcolor=black, 
      citecolor=black, urlcolor=black]{hyperref}

\numberwithin{equation}{section}

\newtheorem{theorem}[equation]{Theorem}

\newtheorem{proposition}[equation]{Proposition}

\newtheorem{remark}[equation]{Remark}

\begin{document}
\title[Emptiness of  homogeneous linear systems with ten general base points]
{Emptiness of homogeneous linear systems with ten general base points}

\author{Ciro Ciliberto}
\address{Dipartimento di Matematica, II Universit\`a di Roma, Italy}
\email{cilibert@axp.mat.uniroma2.it}
\author{Olivia Dumitrescu, Rick Miranda}
\address{Colorado State University, Department of mathematics, College of Natural Sciences,
117 Statistics Building,
Fort Collins, CO 80523 }
\email{rick.miranda@math.colostate.edu, dumitres@math.colostate.edu}
\author{Joaquim Ro\'e}
\address{Departament de Matem\`atiques, Universitat Aut\`onoma de Barcelona, 
Edifici C, Campus de la UAB, 08193 Bellaterra (Cerdanyola del Vall\`es)}
\email{jroe@mat.uab.cat}

\maketitle

\begin{abstract}
In this paper we give a new proof of the fact that for all pairs of positive integers $(d,m)$ with $d/m< 117/37$, the linear system
of plane curves of degree $d$ with ten general base points of multiplicity $m$ is empty.
 \end{abstract}

\section*{Introduction}
We will denote by ${\mathcal{L}}_{d}(m_{1}^{s_{1}},...,m_{n}^{s_{n}})$
the linear system of plane curves of degree $d$
having multiplicities at least $m_{i}$ at $s_{i}$ fixed points, $i=1,\ldots,n$. 
The points in question may be proper or infinitely near, but 
often we will assume them to be general.
In the \emph{homogeneous case},
he \emph{expected dimension} of  the linear system ${\mathcal{L}}_{d}(m^{n})$ is
$$e(\mathcal{L}_{d}(m^{n}))=\max\{-1, \frac{d(d+3)}{2}-\frac{nm(m+1)}{2}\}.$$

Nagata's conjecture for ten general points states that
if $\frac{d}{m}<\sqrt{10} \approx 3.1622$ then ${\mathcal{L}}_{d}(m^{10})$ is empty.
Harbourne and Ro\'e \cite {HR04} proved that
if $\frac{d}{m}<177/56 \approx 3.071$ then ${\mathcal{L}}_{d}(m^{10})$ is empty.
Then Dumnicki \cite{Dum08} (see also \cite{BRHKKSS08}),
combining various techniques, among which
methods developed by Ciliberto--Miranda \cite {CM00} and Harbourne--Ro\'e,
found a better bound $313/99 \approx 3.161616$.
The aim of this paper is to develop a general 
degeneration technique for analysing the emptiness of 
${\mathcal{L}}_{d}(m^{n})$ for general points, and we demonstrate it here
in the case $n=10$. 
This technique is based on the \emph{blow--up and twist method}
introduced in this setting by Ciliberto and Miranda in \cite {CM00}.
Using this, and precisely exploiting a suitable 
degeneration of the plane blown up at ten general points into a union of nine surfaces, 
we prove that ${\mathcal{L}}_{d}(m^{10})$ is empty
if $\frac{d}{m}<\frac{117}{37} \approx 3.162162$.
Using the same degeneration
Ciliberto and Miranda recently proved in \cite {CM08} the non-speciality of ${\mathcal{L}}_{d}(m^{10})$
for $\frac{d}{m} \geq \frac{174}{55}$ and,
as remarked in that article,
one obtains as a consequence the emptyness of ${\mathcal{L}}_{d}(m^{10})$
for $\frac{d}{m} < \frac{550}{174} \approx 3.1609$.
Our emptiness result implies that the $10$--point  Seshadri constant
of the plane is at least $117/370$ (see \cite{HR04}). 
Recently T. Eckl \cite {E09} also obtained the same bound. Using the methods developed in \cite {CM08}
he constructs a more complicated degeneration of the plane into $17$ surfaces 
to find the bound $370/117$
for asymptotic  non--speciality of ${\mathcal{L}}_{d}(m^{10})$.
As proved in \cite  {CM08} this is equivalent to saying that the Seshadri constant 
has to be at least $117/370$,
which is the same conclusion we obtain here with considerably less effort.

The present paper has to be considered as a continuation of \cite {CM08},
which the interested reader is  encouraged to consult 
for details on which we do not dwell here.
From  \cite {CM08} we will take the general setting and most of the notation.
Indeed, the degeneration we use here has been introduced in \cite {CM08}, \S 9.
It is a family parametrized by a disk whose general member 
 $X_t$ is a plane blown up at ten general points, 
whereas the central fibre $X_0$ is 
a local normal crossings union of nine surfaces. 
This construction is briefly reviewed in \S \ref {sec2}. 

A \emph{limit line bundle} on $X_0$ is the datum of a line bundle on the normalization of
each component,
verifying \emph{matching conditions}, i.e. the line bundles have to
agree on the double curves of $X_0$.
In order to analyse the emptiness of ${\mathcal{L}}_{d}(m^{10})$ in the 
asserted range, we use the concept of \emph{central effectivity} introduced in 
 \cite {CM08}, \S 10.1. 
A line bundle $\mathcal L_0$  on $X_0$ is \emph{centrally effective}
if a general section of  $\mathcal L_0$  does not vanish identically on any 
irreducible component of $X_0$.  In particular, if $\mathcal L_0$ is centrally
effective then its restriction to  each component of $X_0$ is effective. 
If ${\mathcal{L}}_{d}(m^{10})$ is not empty, then there is a line bundle
$\mathcal{L} $ on the total space $X$ of the family with a non--zero 
section $s$
vanishing on a surface whose restriction to  the general fiber $X_t$
is a curve in ${\mathcal{L}}_{d}(m^{10})$.  Then  
there is a limit curve in the central fiber $X_0$ as well,
hence there is a limit line bundle $\mathcal{L}_0$ 
associated to that curve. The bundle $\mathcal{L}_0$, which
is the restriction to $X_0$ of $\mathcal L$ twisted by multiples of the
components of $X_0$ where $s$ vanishes, is centrally effective.
In conclusion,  if ${\mathcal{L}}_{d}(m^{10})  \neq \emptyset$ then
there is a limit line bundle which is centrally effective. 
Conversely if for fixed $d$ and $m$ no limit line bundle $\mathcal{L}_0$
is centrally effective, e.g. if its restriction to some component of $X_0$
is not effective, then we conclude that  ${\mathcal{L}}_{d}(m^{10})=\emptyset$.

In this article we will exploit this argument. 
We will describe in \S \ref {sec3} 
limit line bundles $\mathcal{L}_0$
of the line bundle ${\mathcal{L}}_{d}(m^{10})$.
We will see that, in order to apply the central effectivity argument, we can restrict 
our attention to some \emph{extremal}  limit line bundles, and verify central 
effectivity properties only for them.
In \S \ref {sec3} we will  prove that  ${\mathcal{L}}_{d}(m^{10})$ with general base points 
is empty if if $\frac dm < \frac{117}{37}$, by showing that  
none of the extremal limit line bundles verifies the required
central effective properties. 

\section{The degeneration}
\label{sec2}

Consider $X \to \Delta$ the family
obtained by blowing up a point in the central fiber of 
the trivial family over a disc $\Delta\times {\mathbb{P}}^ 2\to \Delta$.
The general fibre $X_t$ for $t\neq 0$ is a ${\mathbb{P}}^ 2$,
and the central fibre $X_0$ is the union of two surfaces
$V\cup Z$, where $V \cong {\mathbb{P}}^ 2$,
$Z \cong {\mathbb{F}}_ {1}$,
and $V$ and $Z$ meet along a rational curve $E$
which is the $(-1)$--curve on $Z$ and a line on $V$ 
(see Figure 1 in \cite {CM08}). 

Choose four general points on $V$ and six general points on $Z$.
Consider these as limits of ten general points in the general fibre $X_t$
and blow them up in the family $X$ (we abuse notation and denote by $X$ also the new family). 
This creates ten exceptional surfaces
whose intersection with each fiber $X_t$ is a $(-1)$--curve,
the exceptional curve for the blow--up of that point.
The general fibre $X_{t}$ of the new family
is a plane blown up at ten general points.
The central fibre $X_{0}$ is the union of $V_1$ a plane blown up at four general points,
and $Z_1$ a plane blown up at seven general points (see Figure 2 in \cite {CM08}). This is the first 
degeneration in \cite {CM08}, \S 3. 

We will briefly recall the notion of a $2$--throw as described in \cite {CM08}, \S 4.2.
Consider a degeneration of surfaces containing two components $V$ and $Z$,
transversely meeting along a double curve $R$.
Let $E$ be a $(-1)$--curve on $V$  intersecting $R$ transversely twice.
Blow it up in the total space.
This creates a ruled surface $T\cong {\mathbb{F}}_{1}$ meeting
$V$ along $E$; the double curve $V \cap T$ is the negative section of $T$.
The surface $Z$ is blown up twice, with two exceptional divisors $G_1$ and $G_2$.
Now blow up $E$ again, creating a double surface $S \cong{\mathbb{F}}_{0}$
in the central fibre meeting
$V$ along $E$
and $T$ along the negative section. 
The blow--up affects $Z$,
by creating two more exceptional divisors $F_1$ and $F_2$
which are $(-1)$ curves, while 
$G_1$ and $G_2$ become $(-2)$--curves.
Blowing $S$ down by the other ruling contracts $E$ on the surface $V$;
$R$ becomes a nodal curve, and $T$ changes into a plane ${\mathbb{P}}^2$
(see Figure 3 in \cite {CM08}). In this process $Z$ becomes non--normal, since
we glue $F_1$ and $F_2$. However, in order to analyse divisors and line bundles on the resulting surface we will always refer to its normalization $Z$. 

On $Z$ we introduced two pairs of infinitely near points $p_i, q_i$,
corresponding to the $(-1)$--cycles $F_i+G_i$ and $F_i$, $i=1,2$. 
Given a linear system $\mathcal L$
on $Z$,  denote by $\mathcal L$ also its pull--back on the blow--up
and consider the linear system $\mathcal L(-a(F_i+G_i)-bF_i)$. We will
say that this system is obtained by imposing to $\mathcal L$ a \emph{point of type} $[a,b]$ at $p_i, q_i$. 

The above discussion is general; we now apply it to the degeneration $V_1\cup Z_1$ described above.
Perform the sequence of $2$--throws along  the following $(-1)$--curves:
\begin{enumerate}
\item The cubic ${\mathcal{L}}_{3}(2,1^6)$ on $Z_1$. This creates the second degeneration in \cite {CM08}, \S 6 (see Figure 5 there).  Note that $V_1$ becomes a $8$--fold blow up of the plane: it started as a $4$--fold blow up and it acquires two more pairs of infinitely near $(-1)$--curves. 
\item Six disjoint curves, i.e.
two conics $C_1={\mathcal{L}}_{2}(1^{4},[1,0],[0,0])$, $C_2={\mathcal{L}}_{2}(1^{4},[0,0],[1,0])$
and four quartics $Q_j=\mathcal{L}_{4}(2^3,1,[1,1]^{2})$ on $V_1$ (the multiplicity one proper point is located at the $i$-th point of the four we blew up on $V$). Trowing the conics creates the third degeneration in \cite {CM08}, \S 7 (see Figure 5 there), and further throwing the quartics creates the fourth degeneration in \cite {CM08}, \S 9 (see Figure 7 there).
\end{enumerate}

By executing all these $2$--throws we introduce 
seven new surfaces
$T$, $U_{i}$, $i=1,2$ (denoted by $T_4$, $U_{i,4}$, $i=1,2$ in \cite {CM08})
and $Y_{j}$, $j=1,\ldots,4$. They are all projective planes, except $T$, which is however a plane at the second degeneration level. Moreover, we have the proper transforms $V$ and $Z$ of $V_1$ and $Z_1$ (denoted $V_4$ and $Z_4$ in \cite {CM08}).  
Throwing  the two conics $C_{i}$
both $Z_1$ and the plane corresponding to $T$ undergo four blow--ups, two of them infinitely near. 
By throwing the four quartics $Q_j$, 
$V_1$ becomes more complicated with $16$ additional blow ups,
in eight pairs of infinitely near points.

\section{The limit line bundles}
\label{sec3}

Next we describe the limit line bundles of $\mathcal L_d(m^ {10})$. Their restrictions to the components of the central fibre will in general be of the form
$$ {\mathcal{L}}_{Z}={\mathcal{L}}_{d_Z}(\mu, q^6, [x_{i},x'_{i}]_{i=1,2}), \quad {\mathcal{L}}_{V}={\mathcal{L}}_{d_V}(\nu^ 4, [y,y']^ 2, [z_i,z'_i]^ 2_{i=1,\ldots,4})$$
$${\mathcal{L}}_{T}={\mathcal{L}}_{d_T}([x_{i},x'_i] _{i=1,2}),\quad {\mathcal{L}}_{U_i}={\mathcal{L}_{s_i}}, i=1,2, \quad 
{\mathcal{L}}_{Y_i}={\mathcal{L}_{t_i}}, i=1,\ldots,4$$
where the parameters $d_Z, \mu, q, x_i, x'_i, ... $ etc. are integers.  
Note that in $\mathcal L_Z$ and $\mathcal L_V$ the points are no longer in general position,
since they have to respect constraints dictated by the 2--throws. 

The matching conditions involving the $U_i$'s and the $Y_i$'s,
imply $s_i=x_i-x'_i$, $i=1,2$,  and  $t_i=z_i-z'_i$, $i=1,\ldots, 4$.  Next we have to impose the remaining matching conditions and also the conditions that this is a limit line bundle of $\mathcal L_d(m^ {10})$, i.e. conditions telling us that the total degree of the limit bundle is $d$ and the multiplicity at the original blown up points is $m$. This would give us the form of all possible limits line bundles of $\mathcal L_d(m^ {10})$, that we need in order to apply the central effectivity argument. However we can simplify our task, by making the following remark.

Let us go back to the 2--throw construction. Let $\mathcal L$ be an effective line bundle 
on the total space of the original degeneration such that $\mathcal L\cdot E=-\sigma<0$. Assume $\sigma=2h$ is even (this will be no restriction in our setting). Create the two exceptional surfaces $S$ and $T$ and still denote by $\mathcal L$ 
the pull--back of the line bundle on the new total space. 
In order to make it centrally effective we have to twist it to $\mathcal L(-uT-(u+v)S)$, and central effectivity requires $u\geq h$, $u\geq v\geq 0$ and $u+v\geq 2h$ (see \cite {CM05}, \S 2). The main remark is that in our setting we may assume $u+v=2h$ by replacing $(u,v)$ with $(u',v')$ where $u'=\min \{u,2h\}$, $v'=2h-u'$. Indeed, $u+v>2h$ means subtracting $E$ more than $2h$ times from $\mathcal L_V$, and creating points of type $[u,v]$ rather than $[u', v']$ for $\mathcal L_Z$. In both cases, this imposes more conditions on the two systems. This is clear for $\mathcal L_V$. As for $\mathcal L_Z$, this follows from $u(F_i+G_i)+vF_i\geq u'(F_i+G_i)+v'F_i$, $i=1,2$.
Therefore if one is able to prove that either one of the two systems on $V$ and $Z$
 is empty, the central effectivity argument will certainly apply to the original twist $\mathcal L(-uT-(u+v)S)$. Note that  $u+v=2h$  is equivalent to require that $\mathcal L(-uT-(u+v)S)\cdot E=0$. 
Essentially the same argument shows that we can also assume that $(u,v)=(h,h)$.

The above discussion shows that, in particular, we may assume $x_{i}=x'_{i}$, ${i=1,2}$, $y=y'$, and  $z_i=z'_i$, ${i=1,\ldots,4}$, with the further conditions that the restrictions to the the 2--thrown curves have degree 0. We call  \emph{extremal}  the bundles verifying these conditions. If, for given $d$ and $m$,  for all extremal limit line bundles either $\mathcal L_Z$ or $\mathcal L_V$ are empty, then there is no centrally effective limit line bundle
and therefore $\mathcal{L}_{d}(m^{10})$ is empty for general points. 

For an extremal bundle, matching between $V$ and $T$ says that $d_T=2x_1=2x_2$. So we set $x_1=x_2=x$. 
The multiplicity conditions for the general points on $V$ then  read
$$m=\nu+4x+2z_i+4\sum_{j\neq i}z_j,\quad  i=1,\ldots,4$$
 yielding $z_1=\ldots=z_4$, which we denote by $z$. Thus we have eight parameters $d_V, d_Z,\nu,\mu, q, x, y, z$ subject to the following seven linear equations
$$3d_Z-2\mu-6q=2d_V-4\nu-y=4d_V-7\nu-4y=0$$
$$m=\nu+4x+14z=q+2x+16z+2y,\quad d=d_Z+6y+48z+6x,  \quad d_V-4y=\mu-4x.$$
The first three come from  the zero restriction conditions to the 2-thrown curves, the next two from the multiplicity $m$ conditions on $V$ and $Z$, the next one from the degree $d$ condition, the last from the matching between $V$ and $Z$. 

Set $\alpha=d-3m$ and $\ell=19m-6d$. By solving the above linear system, we find
$$d_Z=10\alpha-6a,\quad \mu =6\alpha-3a,\quad q=3\alpha-2a, \quad x=5m-\frac 32 d-a$$
$$d_V=9a-18\ell,\quad \nu=4a-8\ell, \quad y= 2a-4\ell, \quad z=\frac \ell 2 .$$
The solutions, as natural, depend on a parameter $a\in \mathbb Z$ (which is the one introduced in the first degeneration in \cite {CM08}).  They are integers since we may assume $d$ and $m$ to be even. 

In conclusion we proved:

\begin{proposition}\label{extremal} In the above degeneration, the extremal limit line bundles $\mathcal L$ of $\mathcal L_d(m^ {10})$ with general base points restrict to the components of the central fibre $X_0$ as follows

$$ {\mathcal{L}}_{Z}={\mathcal{L}}_{10\alpha-6a}(6\alpha-3a, (3\alpha-2a)^6, [5m-\frac 32 d-a,5m-\frac 32 d-a]^ 2)$$
$${\mathcal{L}}_{V}={\mathcal{L}}_{9a-18\ell}((4a-8\ell)^ 4, [2a-4\ell,2a-4\ell]^ 2, [\frac \ell 2 ,\frac \ell 2]^ 8)$$
$${\mathcal{L}}_{T}={\mathcal{L}}_{10m-3 d-2a}( [5m-\frac 32 d-a,5m-\frac 23 d-a]^ 2),\quad {\mathcal{L}}_{U_i}={\mathcal{L}_{0}}, i=1,2, \quad 
{\mathcal{L}}_{Y_i}={\mathcal{L}_{0}}, i=1,\ldots,4.$$

If for all $a\in \mathbb Z$ either $\mathcal L_Z$ or $\mathcal L_V$ is empty,  then no limit line bundle of $\mathcal L_d(m^ {10})$  on $X_0$ is centrally effective, hence $\mathcal L_d(m^ {10})$ is empty. \end{proposition}

\begin{remark}\label  {cremona}{\rm 
As in \cite {CM08}, it is convenient to consider Cremona  equivalent models of the linear systems $\mathcal L_V$ and $\mathcal L_Z$ appearing in Proposition \ref {extremal}.  

The system $\mathcal L_V$ is Cremona equivalent to $\mathcal L_{a-2\ell}( [\frac \ell 2 ,\frac \ell 2]^ 8)$.
The position of the eight infinitely near singular points is special: there are two conics $\Gamma_1,\Gamma_2$ intersecting at four distinct points (the contraction of the four quartics), and each of them contains four of the infinitely near points. The conics $\Gamma_1,\Gamma_2$  are the proper transforms of $F_1,F_2$. For all this, see \cite {CM08}, Lemma 9.1.

The system $\mathcal L_Z$ is Cremona equivalent to$\mathcal L_{76d-240m-3a}((13d-41m-a)^ 6,(\frac {69}2 d-109m-a)^ 4)$. This reduction follows by Lemma 9.2 of \cite {CM08}, but one has to apply a further quadratic transformation based at the three points of multiplicity $\alpha-\ell-a$ of the system there.}\end{remark}

\section{Proof of the theorem}
\label{sec3}

We can now prove our result:

\begin{theorem}\label{main} If $\frac dm<\frac {117}{37}$ then the linear system $\mathcal L_d(m^ {10})$ with ten general base points is empty.
\end{theorem}

\begin{proof} Fix $d,m$ and assume $\mathcal L_d(m^ {10})\neq \emptyset$. According to Proposition \ref {extremal}, there is an integer $a$ such that both $\mathcal L_V$ and $\mathcal L_Z$ are not empty.

Look at the system  $\mathcal L_V$, or rather at its Cremona equivalent form $\mathcal L_{a-2\ell}( [\frac \ell 2 ,\frac \ell 2]^ 8)$ (see Remark \ref {cremona}).  Consider the curve $\Gamma=\Gamma_1+\Gamma_2$, i.e.  the union of the two conics on which the infinitely near base points are located. 
Blow up these base points. By abusing notation we still denote by $\Gamma$  and $\mathcal L_V$ the proper transform of curve and system. Then $\Gamma$ is a 1--connected curve and $\Gamma^ 2=0$. Since $\mathcal L_V$ is effective, one has $\mathcal L_V\cdot \Gamma\geq 0$, i.e. $a\geq 4\ell$. 

Consider then $\mathcal L_Z$, with its Cremona equivalent form $\mathcal L_{76d-240m-3a}((13d-41m-a)^ 6,(\frac {69}2 d-109m-1)^ 4)$. Since this is effective, we have $76d-240m\geq 3a\geq 12\ell$, yielding $\frac dm\geq \frac {117}{37}$. 
\end{proof}

\end{document}